\documentclass[12pt]{article}
\usepackage{amsmath,amssymb}

\def\LC{\mathcal{L}}
\def\TC{\mathcal{T}}

\def\R{\mathbf{R}}

\def\1{\mathbf{1}}

\newcommand{\Pf}{\par\noindent{\em Proof. }}

\def\al{\alpha}
\def\be{\beta}
\def\pa{\partial}
\def\ep{\epsilon}

\newtheorem{prop}{Proposition}[section]
\newtheorem{theorem}{Theorem}[section]

\newtheorem{corollary}{Corollary}
\newtheorem{remark}{Remark}

\begin{document}
\title{GENERALIZED CONTINUOUS-TIME RANDOM WALKS (CTRW), SUBORDINATION BY HITTING TIMES AND FRACTIONAL DYNAMICS}
\author{Vassili N. Kolokoltsov\thanks{Department of Statistics, University of Warwick,
 Coventry CV4 7AL UK,
 and Moscow Institute of Economics, Russia.
  Email: v.kolokoltsov@warwick.ac.uk}}
\maketitle

\begin{abstract}
Functional limit theorem for continuous-time random walks (CTRW) are
found in general case of dependent waiting times and jump sizes that
are also position dependent. The limiting anomalous diffusion is
described in terms of fractional dynamics. Probabilistic
interpretation of generalized fractional evolution is given in terms
of the random time change (subordination) by means of hitting times
processes.
\end{abstract}


\paragraph{Running Head:} Limit distributions for CTRW.

\section{Introduction}\label{sec1}

Suppose $(X_1,T_1),(X_2,T_2),$... is a sequence of i.i.d. pairs of
random variables such that $X_i\in \R^d$,  $T_i\in \R_+$ (jump sizes
and waiting times between the jumps), the distribution of each
$(X_i,T_i)$ being given by a probability measure $\psi (dx \, dt)$
on $\R^d\times \R_+$. Let
\[
N_t=\max \{ n: \sum_{i=1}^n T_i \le t\}.
\]
The process
\begin{equation}
\label{ctrw}
 S_{N_t}= X_1+X_2+...+X_{N_t}
\end{equation}
is called the continuous time random walk (CTRW) arising from
$\psi$. These CTRW were introduced in \cite{MW} and found numerous
applications in physics and economics (see e.g. \cite{Za},
\cite{MS}, \cite{BKS}, \cite{Kor}, \cite{MK} and references
therein). Of particular interest are the situations, where $T_i$
belong to the domain of attraction of a $\be \in (0,1)$-stable law
and $X_i$ belong to the domain of attraction of a $\al \in (0,2)$
-stable law. The limit distributions of appropriately normalized
sums $S_{N_t}$ were first studied in \cite{Kot} in case of
independent $T_i$ and $X_i$ (see also \cite{KKU}). In \cite{BKK2}
the rate of convergence in double array schemes was analyzed and in
\cite{MS} the corresponding functional limit was obtained, which was
shown to be specified by a fractional differential equations. The
importance of the analysis of the case of dependent $T_i$ and $X_i$
was stressed both in \cite{Kot} and \cite{MS}. Here we address this
issue. Moreover we extent the theory to include possible dependence
of $(T_n,X_n)$ on the current position. Our method is quite
different from those used in \cite{Kot}, \cite{KKU}, \cite{MS}. It
is based on the finite difference approximations to continuous-time
operator semigroups and applies the previous results of the author
from \cite{Ko1} on stable-like processes.

It was noted in \cite{MS} that fractional evolution appears from the
subordination of Levy processes by the hitting times of stable Levy
subordinators. Implicitly this idea was present already in
\cite{SZ}. We are going to develop here the general theory of
subordination of Markov processes by the hitting time process
showing that this procedure leads naturally to (generalized)
fractional evolutions. In particular, in spite of the remark from
\cite{MS} that the method from \cite{SZ} (going actually back to
\cite{MW}) "does not identify the limit process" we shall give a
rigorous probabilistic interpretation of the intuitively appealing
(but rather formal) calculations from \cite{SZ}.

In the next Section we demonstrate our approach to the limits of
CTRW by obtaining simple (but nevertheless seemingly new) limit
theorems for position depending random walks with jump sizes from
the domain of attraction of stable laws. In Section 3 these results
will be extended to double scaled random walks, which are needed for
the analysis of CTRW. Section 4 (which is independent of Section 2
and seems to be of independent interest) is devoted to the theory of
subordination by hitting times. In Section 5 we combine the two bits
of the theory from Sections 3 and 4 giving our main results on CTRW.

Let us fix some (rather standard) notations to be used throughout
the paper. For a locally compact space $X$ we denote by $C(X)$ the
Banach space of bounded continuous functions (equipped with the the
sup-norm) and by $C_{\infty}(X)$ its closed subspace consisting of
functions vanishing at infinity. We denote by $(f,\mu)$ the usual
pairing $\int f(x) \mu (dx)$ between functions and measures.
 By a continuous family of
transition probabilities (CFTP) in $X$ we mean as usual a family
$p(x;dy)$ of probability measures on $X$ depending continuously on
$x\in X$, where probability measures are considered in their weak
topology ($\mu_n \to \mu$ as $n\to \infty$ means that $(f,\mu_n) \to
(f,\mu)$ as $n\to \infty$ for any $f\in C(X)$).

For a measure $\mu(dy)$ in $\R^d$ and a positive number $h$ we
denote by $\mu (dy/h)$ the scaled measure defined via its action
 \[
  \int g(z)\mu(dz/h) =\int g(hy) \mu(dy)
  \]
  on functions $g\in C(\R^d)$.

The capital letters $E$ and $P$ are reserved to denote expectation
and probability. The function $\delta (x)$ is the usual Dirac
function (distribution).

\section{Limit theorems for position dependent random walks}\label{sec2}

For a vector $y\in \R^d$ we shall always denote by $\bar y$ its
normalization $\bar y =y/|y|$, where $|y|$ means the usual Euclidean
norm.

Fix an arbitrary $\al \in (0,2)$. Let $S: \R^d\times S^{d-1} \mapsto
\R_+$ be a continuous non-negative function that is symmetric with
respect to the second variable, i.e. $S(x,y)=S(x,-y)$. It defines a
family of $\al$-stable $d$-dimensional symmetric random vectors
(depending on $x\in \R^d$) specified by its characteristic function
$\phi_x$ with
\begin{equation}
\label{stablelaw}
 \ln \phi_x(p)=\int_0^{\infty} \int_{S^{d-1}}
\left(e^{i(p,\xi)}-1-\frac{i(p,\xi)}{1+\xi^2}\right)
  \frac{d |\xi|}{|\xi|^{1+\al}} S(x,\bar \xi) \, d_S\bar \xi,
\end{equation}
where $d_S$ denotes the Lebesgue measure on the sphere $S^{d-1}$. It
is well known that it can be also rewritten in the form

\[
\ln \phi_x(p)=C_{\al}\int_{S^{d-1}} |(p,\bar \xi)|^{\al} S(x,\bar
\xi) \, d_S\bar \xi
\]
with a certain constant $C_{\al}$.

\begin{remark} There are no obstacles for extending our theory to non-symmetric
stable laws. But working with symmetric laws shorten the formulas
essentially.
\end{remark}

\begin{theorem}\label{th1}
Assume
\[
C_1 \le \int_{S^{d-1}} |(\bar p,s)|^{\al} S(x,s) \, d_Ss \le C_2
\]
for all $p$ with some constants $C_1,C_2$ and that $S(x,s)$ has
bounded derivatives with respect to $x$ up to and inclusive order
$q\ge 3$ (if $\al <1$, the assumption $q\ge 2$ is sufficient).
 Then the pseudo-differential operator
\begin{equation}
\label{gen}
 Lf(x)=\ln \phi_x(\frac{1}{i} \frac{\pa}{\pa x})f(x)
 =\int_0^{\infty}
\int_{S^{d-1}} (f(x+y)-f(x))
  \frac{d|y|}{|y|^{1+\al}} S(x,\bar y) \, d_S\bar y
\end{equation}
generates a Feller semigroup $T_t$ in $C_{\infty} (\R^d)$ with the
space $C^{q-1}(\R^d)\cap C_{\infty}(\R^d)$ being its invariant core.
\end{theorem}
This result is proven in \cite{Ko1} and \cite{Ko2}.

\begin{remark} In \cite{Ko1} it is also shown that this semigroup has
a continuous transition density (heat kernel), but we do not need
it.
\end{remark}

Denote by $Z_x(t)$ the Feller process corresponding to the semigroup
$T_t$. We are interested here in discrete approximations to $T_t$
and $Z_x(t)$.

We shall start with the following technical result.

\begin{prop}
\label{convergencetostables}
 Assume that $p(x;dy)$ is a CFTP in ${\R}^d$ from the
normal domain of attraction of the stable law specified by
\eqref{stablelaw}. More precisely assume that for an arbitrary open
$\Omega \in S^{d-1}$ with a boundary of Lebesgue measure zero
\begin{equation}
\label{nor}
 \int_{|y|>n} \int_{\bar y \in \Omega} p(x;dy) \sim \frac{1}{\al
 n^{\al}} \int_{\Omega} S(x,s)\, d_Ss, \quad n\to \infty,
 \end{equation}
 (i.e. the ratio of the two sides of this formula tends to one as
 $n \to \infty$) uniformly in $x$. Assume also that $p(x,\{0\})=0$
 for all $x$.
Then
\begin{equation}
\label{l1}
 \min(1,|y|^2) p(x,dy/h) h^{-\al} \to \min(1,|y|^2)
 \frac{d|y|}{|y|^{\al +1}}
  S(x,\bar y)d_S \bar y, \quad h\to 0,
  \end{equation}
where both sides are finite measures on $\R^d \setminus \{0\}$ and
the convergence is in the weak sense and is uniform in $x\in \R^d$.
If $\al <1$, then also
\[
\min(1,|y|) p(x,dy/h) h^{-\al} \to \min(1,|y|)
 \frac{d|y|}{|y|^{\al +1}} \int_{\Omega}
  S(x,\bar y)d_S \bar y, \quad h\to 0,
  \]
  holds in the same sense.
\end{prop}

\begin{remark} As the limiting measure has a density with respect to
Lebesgue measure, the uniform weak convergence means simply that the
measures of any open or closed set converge uniformly in $x$.
\end{remark}

\Pf By \eqref{nor}
\[
\int_{|z|>A} \int_{\bar z \in \Omega} p(x;dz/h)h^{-\al}
 =\int_{|y|>A/h} \int_{\bar y \in \Omega} p(x;dy)h^{-\al}
  \sim \frac{1}{\al A^{\al}} \int_{\Omega}
  S(x,s)d_Ss
  \]
  as $h\to 0$. Hence
\[
\int_{A<|z|<B} \int_{\bar z \in \Omega} p(x;dz/h)h^{-\al} \to
\int_A^B \frac{d|z|}{|z|^{\al +1}} \int_{\Omega}
  S(x,s)d_Ss.
  \]
  Hence $p(x;dz/h)h^{-\al}$ converges weakly to
$|z|^{-(\al +1)} d|z|
  S(x,z/|z|)d_S(z/|z|)$ on any set separated from the origin.
  It is easy to see that \eqref{l1} follows now from the uniform bound
  \begin{equation}
  \label{bound}
  \int_{|y|<\epsilon} \min(1,|y|^2) p(x,dy/h)h^{-\al} \le C \epsilon
  ^{2-\al}
  \end{equation}
  with a constant $C$. In order to prove \eqref{bound} let us observe
  that
  \[
  \int_{|y|>n} p(x,dy) \le C n^{-\al}
  \]
  with a constant $C$ uniformly for all $x$ and $n>0$ (in fact it
  holds for large enough $n$ by \ref{nor} and is extended to all
  $n$, because all $p(x,dy)$ are probability measures). Hence for an
  arbitrary $\epsilon <1$ one has
  \[
\int_{|y|<\epsilon} \min(1,|y|^2) p(x,dy/h)h^{-\al}
 =\int_{|z|<\epsilon /h} h^2|z|^2) p(x,dy/h)h^{-\al}.
  \]
  Representing this integral as the countable sum of the integrals
  over the regions
  \[
  \epsilon / (2^{k+1}h) <y \le \epsilon / (2^k h),
  \]
  it can be estimated by
  \[
  \sum_{k=0}^{\infty} h^2 \left(\frac{\epsilon} {2^k h}\right)^2 h^{-\al}
  Ch^{\al} 2^{\al (k+1)} \epsilon ^{-\al} =\sum_{k=0}^{\infty}
  C\epsilon ^{2-\al} 2^{\al} 2^{-(2-\al)k}.
  \]
  This yields \eqref{bound}, since the sum on the r.h.s. converges.

 The improvement
concerning the case $\al <1$ is obtained similarly.

 Consider the jump-type Markov process $Z^h(t)$ generated by
 \begin{equation}
 \label{gendis}
 (L_h f)(x) =\frac{1}{h^{\al}}\int (f(x+hy)-f(x))p(x;dy)
\end{equation}
For each $h$ the operator $L_h$ is bounded in $C_{\infty}(\R^d)$ and
hence specifies a Feller semigroup there. The probabilistic
interpretation of $Z^h(t)$ is as follows. Starting at a point $x$
one waits a random $\theta=h^{-\al}$-exponential time $\tau$ (i.e.
distributed according to $P(\tau >t)=\exp (-t\theta)$) and then
jumps to $x+hY$, where $Y$ is distributed according to $p(x;dy)$.
Then the same repeats starting from $x+hY$, etc. In case when $p$
does not depend on $x$
\[
Z^h(t)=h(Y_1+...+Y_{N_t})
\]
is a normalized random walk with the number of jumps $N_t$ being a
Poisson process with parameter $h^{-\al}$, so that $EN_t=th^{-\al}$.
In particular, the number of jumps $n=N_t \sim th^{-\al}$ for small
$h$ so that $Z^h(1)\sim n^{-1/\al}(Y_1+...+Y_n)$.

\begin{theorem}\label{th2}
The semigroup $T_t^h$ generated by $L_h$ converges to the semigroup
$T_t$ generated by $L$. In particular, the corresponding processes
converge in the sense of finite-dimensional marginal distributions.
\end{theorem}

\begin{remark} Everywhere in this paper we work with the convergence of semigroups only.
However by the standard results (see e.g. Theorem 19.25 in
\cite{Kal}) for Feller processes this convergence is equivalent to
the convergence of the distributions of trajectories in an
appropriate Skorokhod space of c\`adl\`ag paths.
\end{remark}

\Pf By \eqref{gendis}
\[
 (L_h f)(x) =\frac{1}{h^{\al}}\int (f(x+z)-f(x))p(x;dz/h),
 \]
 and by Proposition \ref{convergencetostables} this converges to $Lf(x)$ as $h\to 0$ uniformly in $x$
 for $f\in C_{\infty}(\R^d)\cap C^2(\R^d)$.
 By a well known result (see e.g. \cite{Ma})
  the convergence of the generators on the core of the limiting
  semigroup implies the convergence of semigroups.

 The next result concerns the approximations with a non-random
number of jumps. Define the process $S^h_x(t)=S^h_x([t])$ (by the
square bracket the integer part of a real number was denoted) via
\[
S^h_x(0)=x, \quad S^h_x(1)=x+hY_1, \quad ..., \quad
S^h_x(j)=S^h_x(j-1)+hY_j, ...
\]
where each $Y_j$ is distributed according to $p(S_{j-1},dy)$. If
$p(x;dy)$ does not depend on $x$, then
\[
S_x^h(n)=x+h(Y_1+...+Y_n)
 \]
  is just a standard random walk.

  We like to compare the Feller process $Z_x(t)$ on an arbitrary
  fixed time interval $[0,t_0]$ with the discrete approximations
  $S_x^h(t/\tau)$, when the number of jumps $n=t/\tau$ is connected
  with the scaling parameter $h$ by $\tau=h^{\alpha}$.

\begin{theorem}\label{th3}
Under the assumptions of Theorem \ref{th1} and Proposition
\ref{convergencetostables} for any $f\in C_{\infty}(\R^d)$,
$Ef(S_x^h(t/\tau))$ converges to $T_tf(x)$ uniformly on $t\in
[0,t_0]$, as $\tau =h^{\al} \to 0$. In particular, the processes
$S_x^h(t/\tau )$ converge to $Z_x(t)$ in the sense of
finite-dimensional distributions.
\end{theorem}

\Pf It is enough to prove the required convergence for $f\in
C^2(\R^d)\cap C_{\infty}(\R^d)$ only (by Theorem \ref{th1}). Let
such an $f$ be chosen. Denote $f_k(x)=Ef(S^h_x(k))$. Then by the
Markov property $f_k =R_h^kf$, where the operator $R_h$ is defined
via the formula
 \[
 R_hf(x)=\int f(x+hy)p(x;dy).
\]
Clearly each $R_h$ is a positivity preserving contraction on
$C_{\infty}(\R^d)$. On the other hand, the recurrent equation
$f_k=Rf_{k-1}$ can be rewritten as
\begin{equation}
\label{difscheme} \frac{f_k(x)-f_{k-1}(x)}{\tau}
 =h^{-\al} \int (f_{k-1}(x+hy)-f_{k-1}(x)) p(x;dy).
\end{equation}
And this is a discrete time approximation to the equation
\begin{equation}
\label{stableequation}
 \frac{\pa f}{\pa t} =Lf
\end{equation}
on the functions $f\in C^2(\R^d)\cap C_{\infty}(\R^d)$ (and
differentiable in $t$). Since this scheme is well-posed and stable
(as it is solvable uniquely by the contraction $R_h^n$) and the
solution to \eqref{stableequation} is uniquely defined and preserves
the space $C^2(\R^d)\cap C_{\infty}(\R^d)$ (by Theorem \ref{th1}),
it follows by the standard (and easy to prove) general results (see
e.g. \cite{Sam}) that the solutions to the finite-difference
approximation converge to the solution of \eqref{stableequation}.
Theorem is proved.

In case of $p$ not depending on $x$, Theorem \ref{th3} turns to the
known fact on the convergence of random walks with the distribution
of jumps from the domain of normal attraction of a stable law to the
corresponding stable L\'evy motion.

\section{Double-scaled random walks}\label{sec3}

To apply the developed theory to CTRW we shall need a generalization
with multi-scaled walks that we present now.

We are interested in a process in $\R^d\times \R_+$ specified by the
generator

\[
 \LC f(x,u)
 =\int_0^{\infty}
\int_{S^{d-1}} (f(x+y,u)-f(x,u))
  \frac{d|y|}{|y|^{1+\al}} S(x,u, \bar y) \, d_S\bar y
  \]

\begin{equation}
\label{gen2}
  +\int_0^{\infty}(f(x,u+v)-f(x,u))\frac{1}{v^{1+\be}}w(x,u) dv.
\end{equation}

The following result (and its proof) is a straightforward
generalization of Theorem \ref{th1}.

\begin{theorem}\label{th4}
Assume
\[
C_1 \le \int_{S^{d-1}} |(\bar p,s)|^{\al} S(x,u,s) \, d_Ss \le C_2,
\quad
 C_1 \le w(x,u) \le C_2
\]
with some constants $C_1,C_2$ and that $S(x,s)$ and $w(x,u)$ have
bounded derivatives with respect to $x$ and $u$ up to and inclusive
order $q\ge 3$.
 Then the pseudo-differential operator \eqref{gen2}
 generates a Feller semigroup $\TC_t$ in $C_{\infty} (\R^d \times \R_+)$
 (continuous functions up to the boundary) with the space
$(C^{q-1}\cap C_{\infty})(\R^d \times \R_+)$ being its invariant
core and hence a Feller process $(Y,V)(t)$ in $\R^d\times \R_+$.
\end{theorem}

We shall obtain now the corresponding extension of Theorems
\ref{th2}, \ref{th3}.

\begin{theorem}\label{th5}
Assume $p(x,u;dy dv)$ is a CFTP in $\R^d \times \R_+$, which is
symmetric with respect to the reflection  $y\mapsto -y$ and for
which
\[
p(x,u; \{0\}\times \R_+)+p(x,u; \R^d \times \{0\})=0.
 \]
 Assume also that
the projections belong to the domain of normal attraction of stable
laws; more precisely, that uniformly in $(x,u)$
\begin{equation}
\label{nor2}
 \int_{|y|>n} \int_{\bar y \in \Omega} p(x,u;dy dv) \sim \frac{1}{\al
 n^{\al}} \int_{\Omega} S(x,u,s)\, d_Ss, \quad n\to \infty,
 \end{equation}
 and
\begin{equation}
\label{nor3}
 \int_{v>n} \int_{|y|>A} p(x,u;dy dv) \sim \frac{1}{\beta
 n^{\beta}}w(x,u,A), \quad n\to \infty,
 \end{equation}
 for any $A\ge 0$ with a measurable function $w$ of three arguments such that
 \begin{equation}
 \label{limitmeasure}
 w(x,u,0)= w(x,u), \quad \lim_{A\to \infty} w(x,u,A)=0
 \end{equation}
 (so that $w(x,u,A)$ is a measure on $\R_+$ for any $x,u$).

 Consider the jump-type processes generated by
 \begin{equation}
 \label{gendis2}
 (\LC_{\tau} f)(x,u)
  =\frac{1}{\tau}\int (f(x+\tau^{1/\al}y, u+\tau^{1/\be}v))-f(x,u))
  p(x,u;dy dv).
\end{equation}

Then the Feller semigroups $\TC^h_t$ in $C_{\infty}(\R^d\times
\R_+)$ of these processes (which are Feller, because $\LC_h$ is
bounded in $C_{\infty}(\R^d\times \R_+)$ for any $h$) converge to
the semigroup $\TC_t$.
\end{theorem}

\Pf As in Proposition \ref{convergencetostables} one deduces from
\eqref{nor2}, \eqref{nor3} that uniformly in $x,u$
\begin{equation}
\label{asympt1}
 \min(1,|y|^2) \int_0^{\infty} p(x,u;dy/h \, dv) h^{-\al} \to \min(1,|y|^2)
 \frac{d|y|}{|y|^{\al +1}}
  S(x,\bar y)d_S \bar y, \quad h\to 0,
  \end{equation}
and
\begin{equation}
\label{asympt2}
 \min(1,v) \int_{|y|>A} p(x,u;dy dv/h) h^{-\be} \to \min(1,v)
 w(x,u,A) \frac{dv}{v^{\be +1}}, \quad h\to 0,
  \end{equation}

 Next, assuming $f\in (C^2 \cap C_{\infty}(\R^d \times \R_+)$ and writing
\[
\LC_{\tau}f(x,u)=I+II
\]
with
\[
I=\frac{1}{\tau}\int (f(x+\tau^{1/\al}y,u)-f(x,u)) p(x,u;dy dv) +
 \frac{1}{\tau}\int (f(x,u+\tau^{1/\be}v)-f(x,u)) p(x,u;dy dv)
 \]
 and
 \[
 II=\frac{1}{\tau}\int [(f(x+\tau^{1/\al}y,u+\tau^{1/\be}v)-f(x+\tau^{1/\al}y,u))
  -(f(x,u+\tau^{1/\be}v)-f(x,u))] p(x,u;dy dv)
 \]
 one observes that, as in the proof of Theorem \ref{th2},
  \eqref{asympt1} and \eqref{asympt2} (the latter with $A=0$) imply that
 $I$ converges to $\LC f(x,u)$ uniformly in $x,u$. Thus in order to
 complete our proof we have to show that the function $II$ converges to zero, as $\tau \to 0$.
 We have
 \[
 II=\int (g(x+\tau^{1/\al}y,u,v)-g(x,u,v)) p(x,u;dy
 dv/\tau^{1/\be})\frac{1}{\tau}
 \]
 with
 \[
 g(x,u,v)=f(x,u+v)-f(x,u).
 \]
 By our assumptions on $f$
 \[
 |g(x,u,v)| \le C \min (1,v) (\max |\frac {\pa f}{\pa u}|+\max |f|)
  \le \tilde C \min (1,v),
  \]
  and
 \[
 |\frac {\pa g}{\pa x}(x,u,v)| \le C \min (1,v) (\max |\frac {\pa ^2f}{\pa u \pa x}|
 +\max |\frac {\pa f}{\pa x}|)
  \le \tilde C \min (1,v)
  \]
  with some constants $C$ and $\tilde C$. Hence by \eqref{asympt2}
  and \eqref{limitmeasure} for an arbitrary $\ep >0$ there exists a
  $A$ such that
\[
 \int_{|y| >A} (g(x+\tau^{1/\al}y,u,v)-g(x,u,v)) p(x,u;dy dv/\tau^{1/\be})
 \frac{1}{\tau}
 <\ep;
 \]
 and on the other hand, for an arbitrary $A$
 \[
 \int_{|y|<A} (g(x+\tau^{1/\al}y,u,v)-g(x,u,v)) p(x,u;dy dv/\tau^{1/\be})
 \frac{1}{\tau}
 \le \tau^{1/\al}A \kappa
 \]
 with a constant $\kappa$ so that $II$ can be made arbitrary small
 by first choosing large enough $A$ and then choosing small enough
 $\tau$.

 Define now the process $(Y, V)^{\tau}_{x,u}(t/ \tau)
 =(Y, V)^{\tau}_{x,u}([t/ \tau])$, where
\[
(Y, V)^{\tau}_{x,u}(0)=(x,u), \quad
 (Y, V)^{\tau}_{x,u}(1)=(x+\tau^{1/\al}Y_1,u+\tau^{1/\be}V_1),..., \quad
\]
\[
(Y, V)^{\tau}_{x,u}(j)
 =(Y, V)^{\tau}_{x,u}(j-1)+(\tau^{1/\al}Y_j,\tau^{1/\be}V_j), ...
\]
and each pair $(Y_j,V_j)$ is distributed according to $p((Y,
V)^{\tau}_{x,u}(j-1);dy dv)$. If $p(x,u;dy dv)$ does not depend on
$x,u$, then
\[
(Y,V)^{\tau}_{x,u}(n)=(x,u)+(\tau^{1/\al}(Y_1+...+Y_n),\tau^{1/\be}(V_1+...+V_n)).
 \]

In view of Theorem \ref{th5} the following result is obtained by
literally the same arguments as Theorem \ref{th3}.

\begin{theorem}\label{th6}
Under the assumptions of Theorems \ref{th4} and \ref{th5} the linear
contractions $Ef((Y, V)^{\tau}_{x,u}(t/ \tau))$ in $C_{\infty}(\R^d
\times \R_+)$ converge to the semigroup $\TC_tf(x,u)$ of the process
$(Y,V)(t)$ uniformly on $t\in [0,t_0]$, as $\tau \to 0$.
\end{theorem}

\section{Subordination by hitting times and generalized fractional evolutions}\label{sec4}

Let $X(u)$, $u\ge 0$ be a L\'evy subordinator, i.e. an increasing
i.i.d. c\`adl\`ag Feller process (adapted to a filtration on a
suitable probability space) with the generator
\begin{equation}
\label{Levygenerator}
 Af(x)=\int_0^{\infty} (f(x+y)-f(x)) \nu (dy) +a \frac{\pa f}{\pa
 x},
 \end{equation}
 where $a\ge 0$ and $\nu$ is a Borel measure on $\{y>0\}$ such that
 \[
 \int_0^{\infty} \min(1,y) \nu (dy) <\infty.
 \]
 We are interested in the inverse function process or the
 first hitting time process $Z(t)$ defined as
 \begin{equation}
 \label{hittingtimesprocess}
 Z_X(t)=Z(t)=\inf \{ u: X(u) >t\} =\sup \{ u: X(u) \le t\},
 \end{equation}
 which is of course also an increasing c\`adl\`ad process.
 To make our further analysis more transparent (avoiding heavy technicalities of the most general case)
 we shall assume that there exist $\ep >0$ and $\be \in (0,1)$ such that
 \begin{equation}
 \label{conditiononmu}
 \nu (dy) \ge y^{1+\beta}, \quad 0<y<\ep.
 \end{equation}
 For convenient reference we collect in the next statement (without proofs) the
 elementary (well known) properties of $X(u)$.

\begin{prop}
\label{propertiessubordinators} Under condition
\eqref{conditiononmu} (i) the process $X(u)$ is a.s. increasing at
each point, i.e. it is not a constant on any finite time interval;
(ii) distribution of $X(u)$ for $u>0$ has a density $G(u,y)$
vanishing for $y<0$, which is infinitely differentiable in both
variable and satisfies the equation
\begin{equation}
\label{eqondensity} \frac{\pa G}{\pa u}=A^{\star}G,
\end{equation}
where $A^{\star}$ is the dual operator to $A$ given by
\[
A^{\star}f(x)=\int_0^{\infty} (f(x-y)-f(x)) \nu (dy)
  -a \frac{\pa f}{\pa x},
 \]
 (iii) if extended by zero to the half-space $\{t<0\}$ the
 locally integrable function $G(t,y)$ on $\R^2$ specifies a
 generalized function satisfying (in the sense of distribution) the equation
\begin{equation}
\label{eqondensity2} \frac{\pa G}{\pa u}=A^{\star}G+\delta (u)
\delta (y).
\end{equation}
\end{prop}

\begin{corollary}
Under condition \eqref{conditiononmu} (i) the process $Z(t)$ is a.s.
continuous and $Z(0)=0$; (ii) the distribution of $Z(t)$ has a
continuously differentiable probability density function $Q(t,u)$
for $u>0$ given by
\begin{equation}
\label{formulaforQ}
 Q(t,u)=-\frac{\pa}{\pa u} \int_{-\infty}^t
G(u,y) \, dy.
\end{equation}
\end{corollary}

\Pf (i) follows from Proposition \ref{propertiessubordinators} (i)
and for (ii) one observes that
 \[
 P(Z(t) \le u)=P(X(u) \ge t)=\int_t^{\infty} G(u,y) \, dy
 =1- \int_0^t G(u,y) \, dy
 \]
 which implies \eqref{formulaforQ} by the differentiability of $G$.

 \begin{theorem}
 \label{inversesubordinators}
 Under condition
\eqref{conditiononmu} the density $Q$ satisfies the equation
 \begin{equation}
\label{eqonQ1}
 A^{\star}Q=\frac{\pa Q}{\pa u}
 \end{equation}
  for $u>0$,
 where $A^{\star}$ acts on the variable $t$,
 and the boundary condition
 \begin{equation}
 \label{boundarycondition}
 \lim_{u\to 0} Q(t,u)=-A^{\star} \theta (t)
 \end{equation}
where $\theta (t)$ is the indicator function equal one (respectively
$0$) for positive (respectively negative) $t$.
 If $Q$ is extended by
 zero to the half-space $\{u<0\}$, it satisfies the equation
\begin{equation}
\label{eqonQ2}
 A^{\star}Q= \frac{\pa Q}{\pa u}+\delta (u) A^{\star}\theta (t),
\end{equation}
in the sense of distribution (generalized functions).

Moreover the (point-wise) derivative $\frac{\pa Q}{\pa t}$ also
satisfies equation \eqref{eqonQ1} for $u>0$ and satisfies the
equation
\begin{equation}
\label{eqonQ3}
 A^{\star}\frac{\pa Q}{\pa t}= \frac{\pa }{\pa u} \frac{\pa Q}{\pa t}
 +\delta (u) \frac{d}{dt}A^{\star}\theta (t)
\end{equation}
in the sense of distributions.
\end{theorem}

\begin{remark} In the case of a $\be$-stable subordinator $X(u)$
with the generator
\begin{equation}
\label{stablegenerator}
 Af(x)=-\frac{1}{\Gamma (-\be)}\int_0^{\infty} (f(x+y)-f(x))y^{-1-\be} dy,
 \end{equation}
 one has
\begin{equation}
\label{fractional}
 A=-\frac{d^{\be}}{d(-t)^{\be}}, \quad A^{\star}=-\frac{d^{\be}}{dt^{\be}}
\end{equation}
(these equations can be considered as the definitions of fractional
derivatives; we refer to books \cite{MR} and \cite{SW} for a general
background in fractional calculus; a short handy account is given
also in Appendix to \cite{SZ}), in which case equation
\eqref{eqonQ2} takes the
 form
\begin{equation}
\label{eqonQ4}
 \frac{ d^{\be}Q} {dt^{\be}}+ \frac{\pa Q}{\pa u}=\delta (u) \frac{t^{-\be}}{\Gamma (1-\be)}
\end{equation}
coinciding with (B14) from \cite{SZ}.
\end{remark}

\Pf Notice that by \eqref{formulaforQ}, \eqref{eqondensity} and by
the commutativity of the integration and $A^{\star}$ one has
\[
 Q(t,u)=-\int_{-\infty}^t
\frac{\pa}{\pa u} G(u,y) \, dy
 =-\int_{-\infty}^t
(A^{\star} G(u,.))(y) \, dy
  =-A^{\star}\int_{-\infty}^t
 G(u,y) \, dy.
 \]
 This implies \eqref{eqonQ1} (by differentiating with respect to $u$
 and again using \eqref{formulaforQ}) and \eqref{boundarycondition}, because
 $G(0,y)=\delta (y)$.

 Assume now that $Q$ is extended by zero to $\{u<0\}$. Let $\phi$ be
 a test function (infinitely differentiable with a compact support)
 in $\R^2$. Then in the sense of distribution
 \[
 \left((\frac{\pa}{\pa u}-A^{\star})Q,\phi\right)
 =\left( Q, (-\frac{\pa}{\pa u}-A)\phi\right)
 \]
 \[
 =\lim_{\ep \to 0} \int_{\ep}^{\infty} du \int_{\R} dt
  \, Q(t,u)(-\frac{\pa}{\pa u}-A)\phi(t,u)
  \]
  \[
  =\lim_{\ep \to 0} \left[ \int_{\ep}^{\infty} du \int_{\R} dt
  \, \phi(t,u)(\frac{\pa}{\pa u}-A^{\star})Q(t,u)
  +\int_{\R} \phi(t,\ep)Q(t,\ep) \, dt \right].
  \]
  The first term here vanishes by \eqref{eqonQ1}. Hence by \eqref{boundarycondition}
\[
 \left((\frac{\pa}{\pa u}-A^{\star})Q,\phi\right)
 =-\int_{\R} \phi(t,0)A^{\star} \theta (t) \, dt,
 \]
 which clearly implies \eqref{eqonQ2}. The required properties of $\frac{\pa Q}{\pa
 t}$ follows similarly from the representation
 \[
 \frac{\pa Q}{\pa t}(t,u)=-\frac{\pa G}{\pa u}(u,t).
 \]

 We are interested now in the random time change of Markov processes specified by the process
 $Z(t)$.

 \begin{theorem}
 \label{subordinatedprocesses}
 Under the conditions of Theorem \ref{inversesubordinators}
 let $Y(t)$ be a Feller process in $\R^d$, independent of $Z(t)$, and with the domain
 of the generator $L$ containing $(C_{\infty}\cap C^2)(\R^d)$. Denote
 the transition probabilities of $Y(t)$ by
 \[
 T(t,x,dy)=P(Y_x(t) \in dy)=P_x(Y(t)).
 \]

 Then the distributions of the (time changed or subordinated) process $Y(Z(t))$
 for $t>0$ are given by
\begin{equation}
 \label{representationfortransition}
 P_x(Y(Z(t)) \in dy)=\int_0^{\infty} T(u,x,dy)Q(t,u) \, du,
 \end{equation}
 the averages $f(t,x)=Ef(Y_x(Z(t)))$ of $f\in (C_{\infty}\cap C^2)(\R^d)$ satisfy the (generalized)
 fractional evolution equation
\begin{equation}
\label{fracforwardeq}
 A_t^{\star}f(t,x)= -L_xf(t,x)+f(x) A^{\star}\theta (t)
 \end{equation}
(where the subscripts indicate the variables, on which the operators
 act), and their time derivatives $h=\pa f/\pa t$
 satisfy for $t>0$ the equation
\begin{equation}
\label{eqforspeed}
 A_t^{\star}h= -L_xh+f(x) \frac{d}{dt}A^{\star}\theta (t).
 \end{equation}
Moreover, if $Y(t)$ has a smooth
 transition probability density so that $T(t,x,dy)=T(t,x,y)dy$
 and the forward and backward equations
 \begin{equation}
 \label{forback}
 \frac{\pa T}{\pa t}(t,x,y)=L_x T(t,x,y)=L^{\star}_yT(t,x,y)
 \end{equation}
hold, then the distributions of $Y(Z(t))$ have
 smooth density
 \begin{equation}
 \label{representationfordensity}
 g(t,x,y)=\int_0^{\infty} T(u,x,y)Q(t,u) \, du
 \end{equation}
 satisfying the forward (generalized) fractional evolution equation
\begin{equation}
\label{fracforwardeq}
 A^{\star}g= -L_y^{\star}g+\delta (x-y) A^{\star}\theta (t)
 \end{equation}
 and the backward (generalized) fractional evolution equation
\begin{equation}
\label{fracbackwardeq}
 A^{\star}g= -L_xg+\delta (x-y) A^{\star}\theta (t)
 \end{equation}
 with the time derivative $h=\pa g/\pa t$
 satisfying for $t>0$ the equation
\begin{equation}
\label{eqforspeed}
 A^{\star}h= -L_y^{\star}h+\delta (x-y) \frac{d}{dt}A^{\star}\theta (t)
 \end{equation}
\end{theorem}

\begin{remark} In the case of a $\be$-stable L\'evy subordinator $X(u)$
with the generator \eqref{stablegenerator}, where \eqref{fractional}
hold, the left hand sides of the above equations become fractional
derivatives per se. In particular, if $Y(t)$ is a symmetric
$\al$-stable L\'evy motion, equation \eqref{fracforwardeq} takes the
form
\begin{equation}
\label{fracforwardstable}
 \frac{\pa ^{\be}}{\pa t^{\be}}g(t,y-x)
 =\frac{\pa ^{\al}}{\pa |y|^{\al}}g(t,y-x)
 +\delta(y-x)\frac{t^{-\be}}{\Gamma (1-\be)},
\end{equation}
deduced in \cite{SZ} and \cite{U}. The corresponding particular case
of \eqref{representationfordensity} also appears in \cite{MS} as
well as in \cite{SZ}, where it is called a formula of separation of
variables. Our general approach makes it clear that this separation
of variables comes from the independence of $Y(t)$ and the
subordinator $X(u)$ (see Proposition \ref{subordinatedprocesses3}
for a more general situation).
\end{remark}

\Pf For a continuous bounded function $f$ one has for $t>0$ that
\[
Ef(Y_x(Z(t))=\int_0^{\infty} E(f(Y_x(Z(t))|Z(t)=u)Q(t,u) \, du
=\int_0^{\infty} Ef(Y_x(u))Q(t,u) \, du
\]
by the independence of $Z$ and $Y$. This implies
\eqref{representationfortransition} and
\eqref{representationfordensity}.

From Theorem \ref{inversesubordinators} it follows that for $t>0$
\[
A_t^{\star}g=\lim_{\ep \to 0} \int_{\ep}^{\infty}G(u,x,y)
A_t^{\star}Q(t,u) \, du
 =\lim_{\ep \to 0} \int_{\ep}^{\infty}G(u,x,y)
\frac{\pa}{\pa u} Q(t,u) \, du
\]
\[
=- \int_0^{\infty}\frac{\pa}{\pa u} G(u,x,y)Q(t,u) \, du +\delta
(x-y) A^{\star} \theta (t),
\]
where by \eqref{forback} the first term equals $-L^{\star}_yg=L_xg$,
implying \eqref{fracforwardeq} and \eqref{fracbackwardeq}. Other
equations are proved analogously.

Now we like to generalize this theory to the case of L\'evy type
subordinators $X(u)$ specified by the generators of the form
\begin{equation}
\label{Levytypegenerator}
 Af(x)=\int_0^{\infty} (f(x+y)-f(x)) \nu (x,dy) +a (x) \frac{\pa f}{\pa
 x}
 \end{equation}
 with position depending L\'evy measure and drift.
 We need some regularity assumptions in order to have a smooth transition probability density
 like in case of the L\'evy motions.
 \begin{prop}
 \label{Levytypesubordinators}
  Assume that (i) $\nu$ has a density
 $\nu (x,y)$ with respect to Lebesgue measure such that
 \begin{equation}
 \label{conditiononmeasure}
   C_1 \min \left( y^{-1-\be_1}, y^{-1-\be_2}\right) \le \nu(x,y)
    \le C_2 \max \left( y^{-1-\be_1}, y^{-1-\be_2}\right)
\end{equation}
with some constants $C_1,C_2>0$ and
 $0<\beta_1 <\beta_2 <1$ (ii)
  $\nu$ is thrice continuously differentiable with
respect to $x$ with the derivatives satisfying the same estimate
\eqref{conditiononmeasure}, (iii) $a(x)$ is non-negative with
bounded derivatives up to the order three. Then the generator
\eqref{Levytypegenerator} specifies an increasing Feller process
having for $u>0$ a transition probability density $G(u,y)=P(X(u) \in
dy)$ (we assume that $X(u)$ starts at the origin) that is twice
continuously differentiable in $u$.
\end{prop}

\begin{remark} Condition \eqref{conditiononmeasure} holds for
popular stable-like processes with a position dependent stability
index.
\end{remark}

\Pf The existence of the Feller process is proved under much more
general assumptions in \cite{Bass}. A proof of the existence of a
smooth transition probability density is given in \cite{Ko1} under
slightly different assumptions (symmetric multidimensional
stable-like processes), but is easily seen to be valid in the
present situation.

One can see now that the hitting time process defined by
\eqref{hittingtimesprocess} with $X(u)$ from the previous
Proposition is again continuous and has a continuously
differentiable density $Q(t,u)$ for $t>0$ given by
\eqref{formulaforQ}. However \eqref{eqonQ1} does not hold, because
the operators $A$ and integration do not commute. On the other hand,
equation \eqref{eqonQ3} remains true (as easily seen from the
proof). This leads directly to the following partial generalization
of Theorem \ref{subordinatedprocesses}.

 \begin{prop}
 \label{subordinatedprocesses2}
 Let $Y(t)$ be the same Feller process in $\R^d$ as in Theorem
 \ref{subordinatedprocesses2}, but independent hitting time process
  $Z(t)$ be constructed from $X(u)$ under the assumptions of Proposition
  \ref{Levytypesubordinators}.

 Then the distributions of the (time changed or subordinated) process $Y(Z(t))$
 for $t>0$ are given by
 \eqref{representationfortransition} and the time derivatives $h=\pa f/\pa
 t$ of the averages $f(t,x)=Ef(Y_x(Z(t)))$ of continuous bounded functions $f$ satisfy
\eqref{eqforspeed}.
\end{prop}

At last we like to extend this to the case of dependent hitting
times.

\begin{prop}
 \label{subordinatedprocesses3}
 Let $(Y,V)(t)$ be a random process in $\R^d \times \R_+$ such that
(i) the components $Y(t),V(s)$ at different times have a joint
probability density
 \[
 \phi(s,u;y,v)=P(Y(s)\in dy, V(u) \in dv)
 \]
 that is continuously differentiable in $u$ for $u,s>0$, and (ii) the
 component $V(t)$ is increasing and is a.s. not a constant on any
 finite interval. For instance, the process from Theorem \ref{th4}
 enjoys these properties.
  Then (i) the hitting time process $Z(t)=Z_V(t)$ (defined
 by \eqref{hittingtimesprocess} with $V$ instead of $X$) is a.s.
 continuous,
 (ii) there exists a continuous joint probability density of
 $Y(s),Z(t)$ given by
 \begin{equation}
 \label{generaljointdensity}
 g_{Y(s),Z(t)}(y,u)
  =\frac{\pa}{\pa u} \int_t^{\infty} \phi (s,u; y,v) \, dv
  \end{equation}
  and (iii) the distribution of the composition $Y(Z(t))$ has the
  probability density
\begin{equation}
\label{finaldensity}
  \Phi_{Y(Z(t))}(y)
  =\int_0^{\infty} g_{Y(s),Z(t)}(y,s) \, ds
  =\int_0^{\infty}\left(\frac{\pa}{\pa u} \int_t^{\infty} \phi (s,u; y,v) \,
  dv\right)|_{u=s} \, ds.
  \end{equation}
\end{prop}

\Pf (i) and (ii) are straightforward extensions of the Corollary to
Proposition \ref{propertiessubordinators}. Statement (iii) follows
from conditioning and the definition of the joint distribution.

\section{Limit theorems for position dependent CTRW}\label{sec5}

Now everything is ready for our main result.

\begin{theorem}
 \label{mainresult}
Under the assumptions of Theorems \ref{th4} and \ref{th5} let
$Z^{\tau}(t), Z(t)$ be the hitting time processes for
$V^{\tau}(t/\tau)$ and $V(t)$ respectively (defined by the
corresponding formula \eqref{hittingtimesprocess}). Then the
subordinated processes $Y^{\tau}(Z^{\tau}(t)/\tau)$ converge to the
subordinated process $Y(Z(t))$ in the sense of marginal
distributions, i.e.
\begin{equation}
\label{mainconvergence}
 E_{x,0}f(Y^{\tau}(Z^{\tau}(t)/\tau)) \to E_{x,0}(Y(Z(t))),
 \quad \tau \to 0,
 \end{equation}
 for arbitrary $x \in \R^d$, $f \in C_{\infty}(\R^d\times \R_+)$, uniformly for $t$ from any compact interval.
 \end{theorem}

 \begin{remark} We show the convergence in the weakest possible
 sense. It does not seem difficult to extend it to the convergence
 in the
 Skorokhod space of trajectories using standard tools (compactness
 etc) or the theory of continuous compositions from \cite{Wh}.
 Similar result holds for the continuous time approximation from
 Theorem \ref{th5}.
 \end{remark}

 \Pf Since the time is effectively discrete in
 $V^{\tau}(t/\tau)$, it follows that
 \[
 Z^{\tau}(t)=\max \{ u: X(u) \le t\},
\]
and that the events $(Z^{\tau}(t) \le u)$ and $(V^{\tau}(u/\tau) \ge
t)$ coincide, which implies that the convergence of finite
dimensional distributions of $(Y^{\tau}(s/\tau),V^{\tau}(u/\tau))$
to $(Y(s),V(u))$ (proved in Theorem \ref{th6}) is equivalent to the
corresponding convergence of the distributions of
$(Y^{\tau}(s/\tau),Z^{\tau}(t))$ to $(Y(s),Z(t))$.

Next, since $V(0)=0$, is continuous and $V(u)\to \infty$ as $u\to
\infty$ and because the limiting distribution is absolutely
continuous, to show \eqref{mainconvergence} it is sufficient to show
that
\begin{equation}
\label{mainconvergence2}
 P_{x,0}[Y^{\tau}(Z_K^{\tau}(t)/\tau) \in A] \to P_{x,0}[Y(Z_K(t)) \in A], \quad
 \tau \to 0,
\end{equation}
for large enough $K>0$ and any compact set $A$, whose boundary has
Lebesgue measure zero, where
\[
Z_K^{\tau}(t)=Z^{\tau}(t), \quad K^{-1} \le Z^{\tau}(t) \le K,
\]
and vanishes otherwise, and similarly $Z_K(t)$ is defined.

Now
\begin{equation}
\label{proof1}
 P[Y^{\tau}(Z_K^{\tau}(t)/\tau) \in A]=\sum_{k=1/K\tau}^{K/\tau}
 P[V^{\tau}(k)\in A \,\, \& \,\, Z^{\tau}(t)\in [k\tau,(k+1)\tau)]
\end{equation}
and
\begin{equation}
\label{proof2}
 P[Y(Z_K(t)) \in A]=\sum_{k=1/K\tau}^{K/\tau}
 \int_A dy \int_{\tau k}^{\tau(k+1)}g_{Y(s),Z(t)}(y,s) \, ds,
\end{equation}
which can be rewritten as
\begin{equation}
\label{proof3}
 \sum_{k=1/K\tau}^{K/\tau}
 \int_A dy \int_{\tau k}^{\tau(k+1)}g_{Y(\tau k),Z(t)}(y,s) \, ds
 +\sum_{k=1/K\tau}^{K/\tau}
 \int_A dy \int_{\tau k}^{\tau(k+1)}(g_{Y(s),Z(t)}-g_{Y(\tau k),Z(t)})(y,s) \,
 ds.
\end{equation}
The second term here tends to zero as $\tau \to 0$ due to the
continuity of the function \eqref{generaljointdensity}, and the
difference between the first term and \eqref{proof1} tends to zero,
because the distributions of $(Y^{\tau}(s/\tau),Z^{\tau}(t)$
converge to the distribution of $(Y(s),Z(t))$. Hence
\eqref{mainconvergence2} follows. Theorem is proved.

In the case when $S$ does not depend on $u$
 and $w$ does not depend on $x$ in \eqref{gen2}, the limiting
 process $(Y,V)(t)$ has independent components so that the averages
 of the limiting subordinated process satisfy the generalized
 fractional evolution equation from Proposition
 \ref{subordinatedprocesses2},
 and if moreover $w$ is a constant, they satisfy the
 fractional equations from
 Theorem \ref{subordinatedprocesses}.
 In particular, if $p(x,u,dy dv)$ does not depend on $(x,u)$ and
 decomposes into a product $p(dy)q(dv)$, and the limit $V(t)$ is stable, we recover
 the main result from \cite{MS} (in a slightly less general setting,
 since we worked with symmetric stable laws and not with operator
 stable motions as in \cite{MS}), as well as of course the corresponding results
 from \cite{Kot}, \cite{KKU} (put $t=1$ in \eqref{mainconvergence})
 on the long time behavior of the normalized subordinated sums
 \eqref{ctrw}.

{\bf Acknowledgements.} The author is grateful to V. Yu. Korolev, V.
E. Bening and V.V Uchaikin for inspiring him with the beauty of
CTRW, and to J.Hutton and J.Lane for a nice opportunity to deliver a
lecture on CTRW at Gregynog Statistics Workshop 2007.

\end{document}